\documentclass[12pt,fleqn]{article}
\usepackage{graphicx}
\usepackage{pictex}
\usepackage{latexsym}
\renewcommand{\baselinestretch}{1.12}
\newcommand{\R}{{\bf R}}
\newcommand{\C}{{\bf C}}
\newcommand{\Z}{{\bf Z}}
\newcommand{\bP}{{\bf P}}
\newcommand{\be}{{\bf e}}
\newcommand{\cO}{{\cal O}}
\newcommand{\cI}{{\cal I}}
\newcommand{\cM}{{\cal M}}
\newcommand{\cT}{{\cal T}}
\newcommand{\cF}{{\cal F}}
\newcommand{\cE}{{\cal E}}
\newcommand{\cD}{{\cal D}}
\newcommand{\cA}{{\cal A}}
\newcommand{\tN}{\tilde{N}}
\newcommand{\tM}{\tilde{M}}
\newcommand{\tn}{\tilde{n}}
\newcommand{\tm}{\tilde{m}}
\newcommand{\Hom}{\mathop{\rm Hom}\nolimits}
\newcommand{\Spec}{\mathop{\rm Spec}\nolimits}
\newcommand{\vol}{\mathop{\rm vol}\nolimits}
\newcommand{\coker}{\mathop{\rm coker}\nolimits}
\newcommand{\Boxtensor}{\times\hspace{-1em}\Box}
\newtheorem{problem}{Problem}
\newtheorem{proposition}[problem]{Proposition}
\newtheorem{example}[problem]{Example}
\newtheorem{remark}[problem]{Remark}
\title{Problems on \\ Minkowski sums of convex lattice polytopes}
\author{Tadao Oda\\
{\footnotesize {\tt odatadao@math.tohoku.ac.jp}}}
\date{
{\small Abstract submitted at the Oberwolfach Conference\\ ``Combinatorial Convexity
and Algebraic Geometry'' 26.10--01.11, 1997}
}
\begin{document}

\fboxrule=1pt
\thicklines
\maketitle

Throughout, we fix the notation $M:=\Z^r$ and $M_{\R}:=\R^r$.

Given convex lattice polytopes $P,P'\subset M_{\R}$, we have
\[
(M\cap P)+(M\cap P')\subset M\cap(P+P'),
\]
where $P+P'$ is the Minkowski sum of $P$ and $P'$, while the left hand
side means $\{m+m'\mid m\in M\cap P, m'\in M\cap P'\}$.

\begin{problem} \label{problem_1}
\rm
For convex lattice polytopes $P,P'\subset M_{\R}$ when do we have
the equality 
\[
(M\cap P)+(M\cap P')=M\cap(P+P')?
\]
\end{problem}

We always have the equality if $r=1$.
This need not be the case, however, 
if $r\geq 2$ as the following example shows:

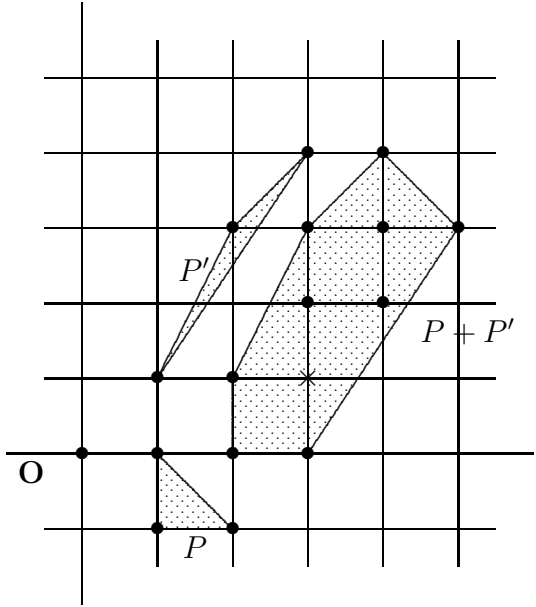
\[
\beginpicture
\setcoordinatesystem units <1cm,1cm>
\setplotarea x from -2 to 7, y from -2 to 6
\putrule from 0 -2 to 0 6
\putrule from 1 -1.5 to 1 5.5
\putrule from 2 -1.5 to 2 5.5
\putrule from 3 -1.5 to 3 5.5
\putrule from 4 -1.5 to 4 5.5
\putrule from 5 -1.5 to 5 5.5
\putrule from -.5 -1 to 5.5 -1
\putrule from -1 0 to 6 0
\putrule from -.5 1 to 5.5 1
\putrule from -.5 2 to 5.5 2
\putrule from -.5 3 to 5.5 3
\putrule from -.5 4 to 5.5 4
\putrule from -.5 5 to 5.5 5
\setlinear \plot 1 -1 2 -1 1 0 1 -1 /
\setlinear \plot 1 1 3 4 2 3 1 1 /
\setlinear \plot 2 0 3 0 5 3 4 4 3 3 2 1 2 0 /
\put {$\bullet$} at 1 -1
\put {$\bullet$} at 2 -1
\put {$\bullet$} at 1 0
\put {$\bullet$} at 1 1
\put {$\bullet$} at 3 4
\put {$\bullet$} at 2 3
\put {$\bullet$} at 2 0
\put {$\bullet$} at 2 1
\put {$\bullet$} at 3 0
\put {$\times$} at 3 1
\put {$\bullet$} at 3 2
\put {$\bullet$} at 3 3
\put {$\bullet$} at 4 2
\put {$\bullet$} at 4 3
\put {$\bullet$} at 4 4
\put {$\bullet$} at 5 3
\put {$\bullet$} at 0 0
\put {{\bf O}} [rt] at -.5 -.1
\put {{$P$}} [t] at 1.5 -1.1
\put {{$P'$}} [rb] at 1.7 2.3
\put {{$P+P'$}} [lt] at 4.5 1.8
\shaderectangleson
\setlinear
\setshadegrid span <2pt>
\vshade 1 -1 0 2 -1 -1 /
\vshade 1 1 1 <,z,,> 2 2.5 3 <z,,,> 3 4 4 /
\vshade 2 0 1 <,z,,> 3 0 3 <z,,,> 4 1.5 4 <,z,,> 5 3 3 /
\shaderectanglesoff
\endpicture
\]

In this example, each of $P$ and $P'$ is nice (known as basic or unimodular),
but their relative position is not.

We may regard the case
$P'=\nu P$ for a positive integer $\nu>0$
as a special case of {\em nice relative position\/}.
We have $P+\nu P=(\nu+1)P$, and
\[
(M\cap P)+(M\cap\nu P)\subset M\cap(\nu+1)P.
\]
\begin{problem} \label{problem_2}
\rm
When do we have the equality
\[
(M\cap P)+(M\cap\nu P)=M\cap(\nu+1)P
\]
for all $\nu\in\Z_{>0}$?
\end{problem}

This problem is related to the projective normality of projective
toric varieties.

We obviously have the equality if $r=1$.
Koelman \cite{koelman} showed that the equality always holds if $r=2$.

More generally, Sturmfels \cite{sturmfels} and others showed that
the equality holds
if $P$ has a basic (also known as unimodular) triangulation.

In view of toric geometry, the following could be a reasonable
formulation of the problem as to when $P$ and $P'$ are in 
nice relative position:
We fix an $r$-dimensional convex lattice polytope $P$, and let
$P'$ to be obtained from $P$ by {\em independent parallel translation\/}
of facets (codimension one faces) of $P$. 
The combinatorial face structure of $P'$ might differ from that
of $P$.

By toric geometry, the convex lattice polytope $P$ corresponds to an 
$r$-dimensional projective toric variety $X$ over the complex number 
field $\C$ together with 
an ample divisor $D$ on $X$,
while $P'$ gives rise to an effective divisor $D'$ on $X$.
$D'$ is ample if the combinatorial face structure of $P'$ coincides 
with that of $P$. When $D'$ is merely nef, the combinatorial face structure
of $P'$ could be slightly degenerate.

\begin{problem} \label{problem_3}
\rm
If $D'$ is nef, do we have the surjectivity of the canonical multiplication map
\[
H^0(X,\cO_X(D))\otimes_{\C} H^0(X,\cO_X(D'))\longrightarrow
H^0(X,\cO_X(D+D'))?
\]
\end{problem}

We know that $H^1(X,\cO_X(D'))=0$ when $D'$ is nef, hence
\[
H^1(X\times X,\cO_{X\times X}(p_1^{-1}D+p_2^{-1}D'))=0
\] 
by K\"unneth formula. 
Consequently, Problem \ref{problem_3} is equivalent to the following:

\begin{problem} \label{problem_4}
\rm
Let $\cI$ be the $\cO_{X\times X}$-ideal corresponding to the
diagonal subvariety $\Delta(X)$ of $X\times X$. 
If $D'$ is nef, do we have 
\[
H^1(X\times X,\cI\otimes_{\cO_{X\times X}}
\cO_{X\times X}(p_1^{-1}D_1+p_2^{-1}D'))=0?
\]
\end{problem}

Hopefully, we might have an affirmative answer at least when
$X$ is {\em smooth\/} and $D'$ is {\em ample}.
There have been unsuccessful attempts in this direction by means of 
Frobenius splittings in characteristic $p>0$.

Without assuming $X$ to be smooth nor $D'$ to be ample, let us give
another formulation for the problem.

Let $N:=\Hom_{\Z}(M,\Z)$ with the canonical bilinear pairing
$\langle\phantom{m},\phantom{m}\rangle:M\times N\longrightarrow\Z$.
Consider the finite complete fan $\Sigma$ for $N$ corresponding to $X$.
As usual, denote by 
\[
\Sigma(1)=\{\rho_1,\rho_2,\ldots,\rho_l\}
\]
the set of one-dimensional cones in $\Sigma$, and let
$n_j\in N$ be the primitive generator for $\rho_j\in\Sigma(1)$.
Let us introduce a free $\Z$-module $\tN:=\bigoplus_{j=1}^l\Z\tn_j$
with the basis consisting of the symbols $\{\tn_1,\ldots,\tn_l\}$
corresponding to $\Sigma(1)$, and the $\Z$-linear map
\[
\pi:\tN\longrightarrow N\qquad\mbox{with }\pi(\tn_j):=n_j
\mbox{ for }j=1,\ldots,l.
\]
Let $\tM:=\Hom_{\Z}(\tN,\Z)$ with the dual basis
$\{\tm_1,\ldots,\tm_l\}$.
Since $\pi$ has finite cokernel, the dual $\Z$-linear map
\[
\pi^{\ast}:M\longrightarrow\tM\qquad\mbox{with }
\pi^{\ast}(m):=\sum_{j=1}^l\langle m,n_j\rangle\tm_j
\mbox{ for any }m\in M
\]
is injective. Let $\cM:=\coker(\pi^{\ast})$ and denote by
$\mu_j\in\cM$ the image of $\tm_j\in\tM$.
We call $(\cM,\{\mu_1,\ldots,\mu_l\})$ the linear Gale transform
of $(N,\{n_1,\ldots,n_l\})$.

$\cM$ is canonically isomorphic to the Weil divisor class group of $X$
(modulo linear equivalence).
Let
\[
\tM\supset\tM_{\geq 0}:=\sum_{j=1}^l\Z_{\geq 0}\tm_j\qquad\mbox{and}\qquad
\cM\supset\cM_{\geq 0}:=\sum_{j=1}^l\Z_{\geq 0}\mu_j.
\]
$\tM_{\geq 0}$ is canonically isomorphic to the semigroup of torus-invariant
effective Weil divisors on $X$.
For $j=1,2,\ldots,l$ we will use $D_j$ and $\tm_j$ interchangeably
to denote the torus-invariant irreducible Weil devisor corresponding
to the one-dimensional cone $\rho_j$.

The homogeneous coordinate ring introduced by Cox, Audin, Delzant, et al.\ 
(cf.\ \cite{cox}) is the semigroup algebra
\[
S:=\C[\tM_{\geq 0}]=\C[x_1,x_2,\ldots,x_l]
\qquad\mbox{with }x_j:=\be(\tm_j)\in S\mbox{ for }j=1,\ldots,l.
\]
We endow the polynomial ring $S$ with the $(\cM_{\geq 0})$-grading
defined by
\[
\deg x_j:=\mu_j\qquad\mbox{for }j=1,\ldots,l.
\]
For $\alpha\in\cM$, we denote by $S_{\alpha}$ the homogeneous part of
degree $\alpha$.

Note that the $(\cM_{\geq 0})$-graded ring $S$ depends only on the
$1$-skeleton $\Sigma(1)$ of $\Sigma$. 
Problem \ref{problem_1} is more or less equivalent to the following:

\begin{problem} \label{problem_5}
\rm
Given $\alpha,\beta\in\cM_{\geq 0}$, when is the multiplication map
\[
S_{\alpha}\otimes_{\C} S_{\beta}\longrightarrow S_{\alpha+\beta}
\]
surjective?
\end{problem}

The fan $\Sigma$ determines the polyhedral cone
\[
C\subset\cM_{\R}:=\cM\otimes_{\Z}\R
\]
spanned by nef divisor classes. The intersection $\cM\cap C^{\circ}$
of $\cM$ with the interior $C^{\circ}$ of $C$ is 
the semigroup of ample divisor classes on $X$.
Then Problems \ref{problem_3} and \ref{problem_4} are almost equivalent
to the following:

\begin{problem} \label{problem_6}
\rm
Is the multiplication map 
$S_{\alpha}\otimes_{\C} S_{\beta}\rightarrow S_{\alpha+\beta}$
surjective if $\alpha\in\cM\cap C^{\circ}$ and $\beta\in\cM\cap C$?
What if $\alpha,\beta\in\cM\cap C^{\circ}$?
\end{problem}

The study of the diagonal ideal sheaf $\cI\subset\cO_{X\times X}$
is important not only in connection with
Problem \ref{problem_4} but in its own right.
In explaining a possible approach to the study, let us follow
the notation of Cox \cite{cox}.

We denote
\[
x^D:=\prod_{j=1}^lx_j^{a_j}\in S\quad\mbox{for }D=\sum_{j=1}^la_j\tm_j
\in \tM_{\geq 0}
\]
and
\[
\deg x^D:=\sum_{j=1}^la_j\mu_j=:[D].
\]
By our convention $D_j=\tm_j$, we have
\[
x^{D_j}=x_j\quad\mbox{and}\quad [D_j]=\mu_j\qquad\mbox{for }j=1,2,\ldots,l.
\]
For each $\alpha\in\cM$ we denote by $\cO_X(\alpha)$ the $\cO_X$-module
corresponding to the degree-shifted graded $S$-module $S(\alpha)$.
We also need the following notation later:
For each $m\in M$ we denote the zero and polar divisors of the character
$\be(m)$ of the torus regarded as a rational function on $X$ by
\[
D^+(m):=
\sum_{\stackrel{\scriptstyle 1\leq j\leq l}{\langle m,n_j\rangle>0}}
\langle m,n_j\rangle D_j \qquad\mbox{and}\qquad
D^-(m):=
\sum_{\stackrel{\scriptstyle 1\leq j\leq l}{\langle m,n_j\rangle<0}}
(-\langle m,n_j\rangle)D_j,
\]
hence $\pi^{\ast}(m)=D^+(m)-D^-(m)$.

We have a canonical homomorphism of 
$(\cM_{\geq 0}\times\cM_{\geq 0})$-graded $\C$-algebras
\[
S\otimes_{\C}S\longrightarrow\C[\cM_{\geq 0}]\otimes_{\C}S
\]
defined by
\[
x^D\otimes x^E\mapsto \be([D])\otimes x^{D+E}\quad
\mbox{for }D,E\in\tM_{\geq 0}.
\]
The ideal $\cI\subset\cO_{X\times X}$ for the diagonal subvariety
$\Delta(X)\subset X\times X$ obviously corresponds to the
$(\cM_{\geq 0}\times\cM_{\geq 0})$-homogeneous ideal
\[
I:=\ker(S\otimes_{\C}S\rightarrow\C[\cM_{\geq 0}]\otimes_{\C}S).
\]
Problems \ref{problem_5} and \ref{problem_6} ask the
surjectivity of the $(\alpha,\beta)$-component 
$S_{\alpha}\otimes_{\C} S_{\beta}\rightarrow 
\be(\alpha)\otimes_{\C} S_{\alpha+\beta}$
of this homomorphism under various conditions on
$\alpha,\beta\in\cM_{\geq 0}$.

$I$ is a homogeneous binomial ideal in the
$(\cM_{\geq 0}\times\cM_{\geq 0})$-graded 
$\C$-algebra
$S\otimes_{\C}S$ depending only on the $1$-skeleton $\Sigma(1)$ of 
the fan $\Sigma$.
We may try to find nice
$(\cM_{\geq 0}\times\cM_{\geq 0})$-graded $S\otimes_{\C}S$-free
resolutions of $I$ to consider Problem \ref{problem_4}.

Identifying the logarithmic derivatives $dx_j/x_j$ with $\tm_j$
for $j=1,2,\ldots,l$ as usual, we get a canonical injective homomorphism
of graded $S$-modules
\[
\Omega_S^1\longrightarrow S\otimes_{\Z}\tM,\qquad
dx_j\mapsto x_j\otimes\tm_j\quad\mbox{for }j=1,2,\ldots,l.
\]
Denote by
\[
\Omega:=\ker(\Omega_S^1\rightarrow S\otimes_{\Z}\tM\rightarrow
S\otimes_{\Z}\cM)
\]
the kernel of the composite of this homomorphism with the canonical
projection $S\otimes_{\Z}\tM\rightarrow S\otimes_{\Z}\cM$.
It is not hard to show that the sheaf $\Omega_X^1$ of Zariski
differential $1$-forms (resp.\ $\bigoplus_{j=1}^l\cO_X(-D_j)$) 
is the $\cO_X$-module associated to the graded $S$-module
$\Omega$ (resp.\ $\Omega_S^1$).
In this way, we get the following well-known result:

\begin{proposition} \label{prop_euler}
{\rm (Generalized Euler exact sequence. cf.\ 
Batyrev-Cox \cite{batyrev-cox})}
\newline
We have an exact sequence of $\cO_X$-modules
\[
0\rightarrow\Omega_X^1\rightarrow\bigoplus_{j=1}^l\cO_X(-D_j)
\rightarrow\cO_X\otimes_{\Z}\cM\rightarrow 0.
\]
\end{proposition}

\begin{remark} \label{rem_arrangement1}
\rm
The graded $S$-module $\Omega$ is generated over $S$ by
\[
\left\{\left.x^Ddx^E-x^Edx^D\;\right|\;D,E\in\tM_{\geq 0},\;
D\sim E\right\},
\]
hence by
\[
\left\{\left.x^{D^+(m)}dx^{D^-(m)}-x^{D^-(m)}dx^{D^+(m)}
\;\right|\; m\in M
\right\}.
\]
The vectors $n_1,n_2,\ldots,n_l\in N$ give rise to an arrangement
$\cA$ of hyperplanes $\{n_j\}^{\perp}\subset M_{\R}$.
A {\em chamber\/} $\Gamma$ for $\cA$ is one of the top-dimensional
polyhedral cones appearing in the partition of $M_{\R}$ induced by
the arrangement $\cA$. If we choose for each chamber $\Gamma$
a set $\Xi_{\Gamma}$ of generators of the semigroup $M\cap\Gamma$,
then $\Omega$ is generated over $S$ by
\[
\left\{x^{D^+(m)}dx^{D^-(m)}-x^{D^-(m)}dx^{D^+(m)}
\;\left|\; m\in\bigcup_{\mbox{\scriptsize chambers }\Gamma}\Xi_{\Gamma}
\right.\right\}.
\]
\end{remark}

\begin{remark} \label{rem_arrangement2}
\rm
As in Remark \ref{rem_arrangement1}, we see that the
``diagonal'' ideal $I\subset S\otimes_{\C}S$ is generated over
$S\otimes_{\C}S$ by
\[
\left\{\left.x^D\otimes x^E-x^E\otimes x^D\;\right|\;D,E\in\tM_{\geq 0},\;
D\sim E\right\},
\]
hence by
\[
\left\{x^{D^+(m)}\otimes x^{D^-(m)}-x^{D^-(m)}\otimes x^{D^+(m)}
\;\left|\; m\in\bigcup_{\mbox{\scriptsize chambers }\Gamma}\Xi_{\Gamma}
\right.\right\}.
\]
\end{remark}

\begin{example} \label{example_projectivespace}
\rm
When $X=\bP^r$ is the projective space, the corresponding polytope
$P$ is a unimodular simplex, and $l=r+1$.
We have $\cM=\Z$, and 
$S_{\alpha}\otimes S_{\beta}\rightarrow S_{\alpha+\beta}$ in
Problem \ref{problem_5} is obviously
surjective for all $\alpha,\beta\in\cM_{\geq 0}$.
Nevertheless, the description of the diagonal ideal $\cI\subset\cO_{X\times X}$
is nontrivial.
By Beilinson \cite{beilinson} we have an exact sequence
\[
0\rightarrow\cO_X(-r)\Boxtensor\Omega_X^r(r)\rightarrow\cdots\rightarrow
\cO_X(-j)\Boxtensor\Omega_X^j(j)\rightarrow\cdots\rightarrow
\cO_X(-1)\Boxtensor\Omega_X^1(1)\rightarrow\cI\rightarrow 0,
\]
where $\cF\Boxtensor\cF':=(p_1^{\ast}\cF)\otimes_{\cO_{X\times X}}
(p_2^{\ast}\cF')$ is the external tensor product on $X\times X$
of $\cO_X$-modules $\cF$ and $\cF'$.
(Thanks are due to Miles Reid for pointing out this result to the author.)
One way of proving this is to note that the $S\otimes_{\C}S$-module
homomorphism
\[
S\otimes_{\C}\Omega_S^1\longrightarrow S\otimes_{\C}S,\qquad
1\otimes dx_j\mapsto x_j\otimes 1\quad\mbox{for }j=1,2,\ldots,r+1
\]
induces a surjection $S\otimes_{\C}\Omega\rightarrow I$.
Thus the Koszul complex arising out of the 
corresponding $\cO_{X\times X}$-homomorphism
(cf., e.g., Eisenbud \cite{eisenbud})
\[
\cO_X(-1)\Boxtensor\Omega_X^1(1)\longrightarrow\cO_{X\otimes X},
\]
whose cokernel is $\cO_{\Delta(X)}$, gives the exact sequence above.

Similarly, the Koszul complex arising out of 
the $S$-homomorphism $\Omega_S^1\rightarrow S$ which sends
$dx_j$ to $x_j$ for $j=1,2,\ldots,r+1$ is nothing but the complex
$(\Omega_S^{\textstyle \cdot},d)$. 
By Remark \ref{rem_arrangement1}, we see that $\Omega$
is the image of $d:\Omega_S^2\rightarrow\Omega_S^1$.
Consequently, we get an exact sequence of $\cO_X$-modules
\[
0\rightarrow
\cO_X(-r-1)^{\oplus{}_{r+1}C_{r+1}}
\rightarrow\cdots\rightarrow
\cO_X(-j)^{\oplus{}_{r+1}C_j}
\rightarrow\cdots\rightarrow
\cO_X(-2)^{\oplus{}_{r+1}C_2}
\rightarrow\Omega_X^1\rightarrow 0,
\]
where 
\[
{}_{r+1}C_j:=\left(\begin{array}{c}r+1\\j\end{array}\right)
\]
are the binomial coefficients.

Alternatively, we could use the Eagon-Northcott complex 
(cf.\ Eagon-Northcott \cite{eagon-northcott}, Kirby \cite{kirby}
and Eisenbud \cite{eisenbud})
for the $2\times 2$-minors of
\[
\left(
\begin{array}{cccc}
x_1\otimes 1&x_2\otimes 1&\ldots&x_{r+1}\otimes 1\\
1\otimes x_1&1\otimes x_2&\ldots&1\otimes x_{r+1}
\end{array}
\right)
\]
to get an exact sequence
\[
0\rightarrow\cE_r\rightarrow\cdots\rightarrow\cE_p\rightarrow\cdots
\rightarrow\cE_1\rightarrow\cI\rightarrow 0,
\]
where 
\[
\cE_p:=\bigoplus_{\stackrel{\scriptstyle j+k=p+1}{j\geq 1,k\geq 1}}
\bigwedge^{p+1}\tM\otimes_{\Z}\cO_{X\times X}(-j,-k)
\]
with
\[
\cO_{X\times X}(\alpha,\beta):=\cO_X(\alpha)\Boxtensor\cO_X(\beta)
\mbox{ for $\alpha,\beta\in\cM=\Z$}.
\]
\end{example}


\end{document}